%% file: paper.tex
\documentclass[10pt,a4paper,twocolumn]{article}

\usepackage[left=1.5cm,right=1.5cm,top=2.5cm]{geometry}
\usepackage{mathtools}

\usepackage[utf8]{inputenc}
\usepackage[T1]{fontenc}
\usepackage{microtype}
\usepackage{lmodern}

\usepackage{hyperref}
\usepackage{orcidlink}

\usepackage{csquotes}

\usepackage[charsperline=57]{jlcode}

\usepackage{authblk}

\usepackage{pgfplots,pgfplotstable}
\usepgfplotslibrary{groupplots,colorbrewer}

\pgfplotsset{
    compat=1.17,
    table/col sep=comma,
    table/search path={data},
    every axis plot/.append style={thick},
    legend style={font=\scriptsize},
    ticklabel style={font=\scriptsize},
    cycle list/Set1-5,
    cycle multiindex* list={
        mark list*\nextlist
        Set1-5\nextlist
        linestyles*\nextlist
    },
    discard if not/.style 2 args={
        x filter/.code={
            \edef\tempa{\thisrow{#1}}
            \edef\tempb{#2}
            \ifx\tempa\tempb
            \else
                
            \fi
        }
    },
    legend cell align={left},
}

\newlength{\mytemp}
\newlength{\mytempb}
\makeatletter
\define@key{lst}{xgobble}{%
  \setlength{\mytemp}{\lst@xleftmargin}
  \addtolength{\mytemp}{-#1\bfem}
  \def\lst@xleftmargin{\mytemp}%
  \setlength{\mytempb}{\lst@framexleftmargin}
  \addtolength{\mytempb}{-#1\bfem}
  \def\lst@framexleftmargin{\mytempb}%
}
\makeatother

\def\pkg#1{\texttt{#1}}
\def\alg#1{\texttt{#1}}

\def\keywordname{{\bfseries \emph Keywords}}%
\def\keywords#1{\par\addvspace\medskipamount{\rightskip=0pt plus1cm
\def\and{\ifhmode\unskip\nobreak\fi\ $\cdot$
}\noindent\keywordname\enspace\ignorespaces#1\par}}

\title{DelayDiffEq: Generating Delay Differential Equation Solvers via Recursive Embedding of Ordinary Differential Equation Solvers}

\author[1]{David Widmann~\orcidlink{0000-0001-9282-053X}}
\author[234]{Chris Rackauckas~\orcidlink{0000-0001-5850-0663}}
\affil[1]{Uppsala University, Sweden}
\affil[2]{Massachusetts Institute of Technology, USA}
\affil[3]{Julia Computing}
\affil[4]{Pumas-AI}

\hypersetup{
    pdftitle = {DelayDiffEq: Generating Delay Differential Equation Solvers via Recursive Embedding of Ordinary Differential Equation Solvers},
    pdfauthor = {David Widmann, Chris Rackauckas},
    pdfkeywords = {Julia, Delay Differential Equation, Scientific Computing},
}

\begin{document}

\twocolumn[
\begin{@twocolumnfalse}

\maketitle

\begin{abstract}

Traditional solvers for delay differential equations (DDEs) are designed around only a single method and do not effectively use the infrastructure of their more-developed ordinary differential equation (ODE) counterparts. In this work we present \pkg{DelayDiffEq}, a Julia package for numerically solving delay differential equations (DDEs) which leverages the multitude of numerical algorithms in \pkg{OrdinaryDiffEq} for solving both stiff and non-stiff ODEs, and manages to solve challenging stiff DDEs. We describe how compiling the ODE integrator within itself, and accounting for discontinuity propagation, leads to a design that is effective for DDEs while using all of the ODE internals. We highlight some difficulties that a  numerical DDE solver has to address, and explain how \pkg{DelayDiffEq} deals with these problems. We show how \pkg{DelayDiffEq} is able to solve difficult equations, how its stiff DDE solvers give efficiency on problems with time-scale separation, and how the design allows for generality and flexibility in usage such as being repurposed for generating solvers for stochastic delay differential equations.

\end{abstract}

\keywords{Julia, Delay Differential Equation, Scientific Computing}

\vspace{\baselineskip}

\end{@twocolumnfalse}
]

\section{Introduction}

In nature many changes do not occur instantaneously, as prominent
examples such as gestation times and incubation periods indicate. This
lead to the rise of delay differential equations (DDEs) in many areas
of mathematical modeling, such as population dynamics, pharmacokinetics,
and control problems~\cite{baker95_issues_numer_solut_evolut_delay_differ_equat}. In contrast
to ordinary differential equations (ODEs), DDEs allow the rate of change
to depend on past conditions of the considered system.

However, while DDEs are a small conceptual change to ODEs, there is a larger
mathematical difference that arises. Specifically, many DDEs have a discontinuity
in the initial condition which, unlike ODEs, will continue to be in the system and
propagate as we describe further in Section~\ref{sec:disc}. Handling these propagations
is crucial to performance and accuracy, and thus the software of this field, such as the
MATLAB \pkg{dde23}~\cite{shampine01_solving_ddes_matlab} and the classic \pkg{RADAR5}~\cite{guglielmi01_implementing_radau}, are completely separate software from those of
ODE integrators. This gives rise to many development and maintenance issues: while in
theory an advance in linear solvers (Newton-Krylov methods, preconditioners, etc.) can
be used in DDE solvers, in practice none of these solvers support these features while
their more-developed ODE counterparts do. This raises the question as to whether DDE
software could be designed to more effectively make use of the ODE infrastructure so that
all developments of the ODE solvers are inherited by the DDE solvers.

To answer this challenge we developed \pkg{DelayDiffEq}, a DDE solver package in Julia~\cite{bezanson2017julia} which uses the ODE solvers of \pkg{OrdinaryDiffEq}~\cite{DifferentialEquations.jl-2017} to develop the DDE solvers. We describe a form
of ``compilation trick'', compiling an ODE integrator inside of an ODE integrator,
which is used to generate the DDE solvers from the ODE solvers. This alone however does
not guarantee the accurate handling of discontinuities, and thus \pkg{DelayDiffEq} then modifies
the stepping to account for these issues. This gives a DDE solver where all methods of the
ODE solver, ones for non-stiff and stiff equations, are directly inherited by the DDE solver
with all of the available options. We show that this generated solver is capable of handling 
known difficult equations to verify its accuracy.%
\footnote{The source code of the paper and its examples are available on \href{https://github.com/SciML/DelayDiffEq.jl/tree/paper}{Github}.}

\section{Background}

\subsection{Definition}

Loosely speaking, a delay differential equation (DDE) is an ordinary
differential equation (ODE) in which the derivative depends on past
values of the state. More concretely, let $t \in \mathbb{R}$ be the 
independent variable, which usually denotes
time, and $x(t) \in \mathbb{R}^n$ be the dependent state at time $t$. Then
the differential equation
\begin{equation}\label{eq:dde}
    x'(t) = f(t, x(t), x(t - \tau)),
\end{equation}
where $f \colon \mathbb{R} \times \mathbb{R}^n \times \mathbb{R}^n$ and
$\tau > 0$, is called a DDE with constant delay $\tau$. In general, 
a delay $\tau$ may also depend on $t$ or even on both $t$ and $x(t)$.
In this case we call the DDE time-dependent or state-dependent, respectively.
Additionally, DDEs may contain multiple delays. A further generalization is
possible by allowing the state to depend also on its derivative at previous
time points, yielding so-called neutral DDEs.

\subsection{History function}

For ODEs, so-called initial value problems (IVPs) can be formulated by specifying
an initial condition $x(t) = x_0$ at an initial time point $t_0$. For DDEs, 
usually declaring an initial value only at time point $t_0$ is not sufficient.
Instead typically a so-called history function, or initial function, has to be
provided for some times $t \leq t_0$. Thus commonly an IVP formulation for the DDE
in eq.~\eqref{eq:dde} is given by
\begin{equation}\label{eq:dde_ivp}
    \begin{aligned}
        x'(t) &= f(t, x(t), x(t - \tau)), \qquad & t &\geq t_0, \\
        x(t) &= x_0(t), \qquad &t &\leq t_0.
    \end{aligned}
\end{equation}

\subsection{Propagated discontinuities}\label{sec:disc}

Even if $x_0$, $f$, and $\tau(t, x(t))$ are smooth, typically the solution $x$
is non-smooth since it has a derivative jump discontinuity at the initial time
point $t_0$, i.e.,
\begin{equation*}
    \lim_{t \to t_0} x'_0(t) \neq \lim_{t \to t_0} x'(t).
\end{equation*}
If there exist time points $t > t_0$ such that the deviated
argument $t - \tau(t, x(t)) < t_0$, then the discontinuity at time point $t_0$
will propagate in time and induce discontinuities at subsequent time points.
For non-neutral DDEs these propagated discontinuities are usually discontinuities
of higher-order derivatives, leading to a so-called smoothing of the solution,
whereas in general for neutral DDEs smoothing does not occur~\cite{bellen03_numer_methods_delay_differ_equat,neves76_charact_jump_discont_state_depend}.

\subsection{Example: Dynamical structure}

\begin{figure}
    \centering
    \input{figures/mackey_glass.tex}
    \caption{Solution of problem~\eqref{eq:mackey_glass} and its time delay embedding.}
    \label{fig:mackey_glass}
\end{figure}
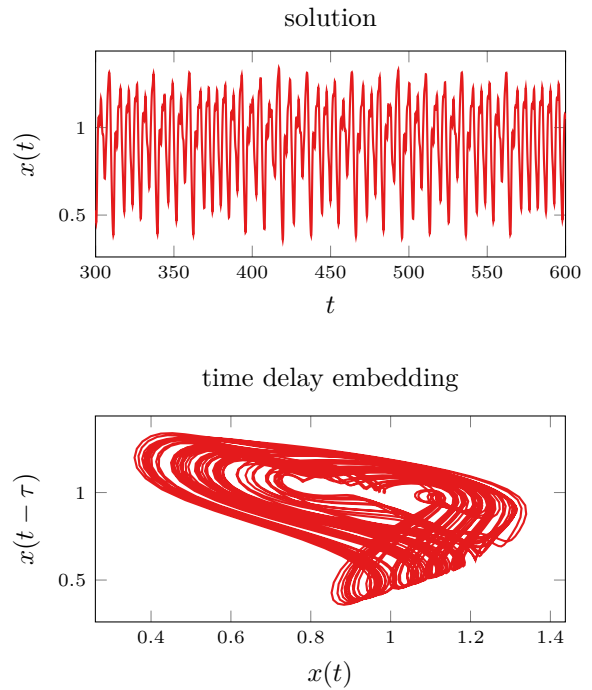

Mackey and Glass' model of circulating blood cells~\cite{glass79} is a scalar DDE system but
still for certain parameter values we can observe a chaotic behavior.
The model is given by
\begin{equation}\label{eq:mackey_glass}
    \begin{aligned}
        x'(t) &= \frac{\lambda \theta^{n} x(t - \tau)}{\theta^{n} + x{(t - \tau)}^n} - \gamma x(t), \qquad & t &\geq 0, \\
        x(t) &= x_{0}, \qquad & t &\leq 0.
    \end{aligned}
\end{equation}
For the parameter values $\lambda = 2$, $\theta = 1$, $n = 9.65$, $\tau = 2$, $\gamma = 1$, and $x_{0} = 0.5$ the system exhibits a chaotic solution~\cite[Figure~6(e)]{glass79}.
In Figure~\ref{fig:mackey_glass} the chaotic solution computed with \pkg{DelayDiffEq} and its time delay embedding are shown for time $300 \leq t \leq 600$.

\subsection{Method of steps}

Consider an IVP for a DDE with one constant delay, given by eq.~\eqref{eq:dde_ivp}.
Assuming that a solution to this problem exists, a solution can be found by
computing solutions of an iterative sequence of IVPs for different ODEs. More
concretely, for $k = 1,2,\ldots$ one computes a solution
$x_k \colon [t_{k-1}, t_{k-1} + \tau] \to \mathbb{R}^n$ of the IVP
\begin{equation}\label{eq:ode_method_of_steps}
    \begin{aligned}
        x'_k(t) &= f(t, x_k(t), x_{k-1}(t - \tau)), \quad t \in [t_{k-1}, t_k], \\
        x_k(t_{k-1}) &= x_{k-1}(t_{k-1}),
    \end{aligned}
\end{equation}
where $t_k \coloneqq t_0 + k\tau$. Then the function $x$ defined by
$x(t) \coloneqq x_k(t) \Longleftrightarrow t \in [t_{k-1}, t_k)$ is a solution
to the original IVP in eq.~\eqref{eq:dde_ivp}.

This algorithm for solving DDEs is called method of steps. It can be applied more
generally also to IVPs of DDEs with multiple state-dependent delays, as long as the
delays are uniformly strictly greater than 0. The same iterative approach can be
used to transfer results from the theory of ODEs such as existence, uniqueness, and
boundedness of solutions to DDEs without vanishing delays.

\section{Numerical challenges}

For DDEs, the computation of numerical solutions of IVPs can be challenging. In
theory the method of steps allows to reduce the problem to the computation of
a sequence of numerical solutions to IVPs of ODEs. A crucial requirement
for this approach to work are so-called continuous ODE methods that provide a
numerical solution not only for a grid of discrete time points but for all
time points in the integration interval.

A disadvantage of the method of steps that it can be inefficient.
The algorithm requires that the step size is upper bounded by the
minimally occurring delay, and hence the number of time steps can become
large if the maximally allowed step size is small relative to the integration time
interval. Additionally, problems with vanishing delays cannot be solved by the
method of steps.

Another difficulty is posed by the fact that usually numerical methods for
solving ODEs expect the solution to be sufficiently smooth in each step. To
ensure this, propagated discontinuities have to be considered and included
in the grid of time points that the numerical method steps to~\cite{bellen03_numer_methods_delay_differ_equat}.
For constant delays the
relevant discontinuities can be computed a priori but for state-dependent
delays this is not possible. Even for constant delays numerical
difficulties can arise, as Shampine and
Thompson~\cite{shampine01_solving_ddes_matlab} illustrate: two constant delays
of $1/3$ and $1$, albeit both inducing propagated discontinuities at time point
$t_0 + 1$, might seem to cause discontinuities at different time points
numerically since $1/3$ cannot be represented exactly by the floating point
arithmetic of computers.

\section{How DelayDiffEq works}

So far one of the core design principles of \pkg{DelayDiffEq} has been to exploit
the large amount of numerical ODE algorithms readily available in the Julia
package \pkg{OrdinaryDiffEq}. The ODE solvers in \pkg{OrdinaryDiffEq} support both stiff
and non-stiff ODEs\footnote{Although it
\textcquote{burrage95_parallel_sequential_methods}{is extremely difficult
to classify stiffness in a completely rigorous manner}, it is commonly characterized
by the existence of processes with significantly different time
scales~\cite{burrage95_parallel_sequential_methods,dahlquist75_recent_work_stiff}.}.
Moreover, the numerical methods are continuous and come with specialized interpolations.

The design of \pkg{DelayDiffEq} is as follows. In order for a user to describe a DDE, we need
to have a history function which the user can access. This history function is a continuous
function of the states of the DDE solution, which, if we only had an ODE, would be given
by the dense output (continuous interpolation) of the ODE solver. Thus what we do is the
following:

\begin{enumerate}
    \item Create an ODE integrator from \pkg{OrdinaryDiffEq} with state size matching the DDE we wish to solve. The integrator comes with the property that $\operatorname{integ}(t)$ computes the interpolation at time point $t$ given the current solution of the ODE.
    \item Generate an ODE by enclosing this interpolation object into the user-defined DDE: $f_{\text{ODE}}(t,x(t)) = f(t,x(t),\operatorname{integ})$.
    \item Create a new ODE integrator to solve the $f_{\text{ODE}}$ system. If the interpolation $\operatorname{integ}$ is correct, then the solution of this new ODE will give the solution of the DDE.
\end{enumerate}

The second ODE integrator will have the correct state values $x(t)$ and stage calculation if the function $\operatorname{integ}$ is correct.
On the other hand, $\operatorname{integ}$ is only a correct interpolation if it has the correct state values $x(t)$ and stage calculations.
Thus we effectively alias the caches between the two versions of the ODE integrators, so that a step of the second ODE integrator updates its history function, which then allows for an accurate step of the integrator.

Notably, this scheme does not require implementing any new numerical integrators and simply relies on the ODE integration itself.
However, as described this scheme is correct only if a new step uses only the history of previous steps, i.e., only when the step size is smaller than the minimum delay time.
Thus while this generates an ODE solver, more work must be done to make this routine efficient.

\subsection{Unconstrained time stepping}\label{sec:unconstrained}

As discussed above, it is desirable to not upper bound the step size that
the numerical solver is allowed to take. However, unconstrained time stepping
comes at a cost. Allowing unconstrained steps implies that the numerical
solution of the DDE is only available as the solution to an implicit problem,
even for explicit ODE methods. Essentially, one must know that $\operatorname{integ}$ is correct to know that the next step is correct, but if
the time step is larger than the minimum delay then the integrator will use an inaccurate extrapolation from the previous time step.
In order to correct for this behavior, we transform this problem into an implicit problem which can be solved via fixed-point iteration.

Formally, in the $k$th integration step we want to compute an approximation
$\hat{x}_k \colon [t_{k - 1}, t_k] \to \mathbb{R}^n$ of the solution to the DDE
in eq.~\eqref{eq:dde} with approximate solution $\hat{x}(t)$ for $t \leq t_{k - 1}$ as history function.
In \pkg{DelayDiffEq}, we perform a fixed-point iteration to compute $\hat{x}_k$ if it is only available implicitly.
Hereby we proceed in the following way:

\begin{enumerate}
    \item \label{enum:extrapolation} Extend the dense solution $\hat{x}$ up to time $t_{k-1}$
    to the time interval $[t_{k-1}, t_k]$ by extrapolating the specialized
    collocation polynomials available in \pkg{OrdinaryDiffEq}.
    \item \label{enum:ode_step} Based on the approximation up to time $t_k$,
    compute a discrete approximation of $x(t_k)$ and additional
    collocation stages in the time interval $[t_{k - 1}, t_k]$ by using an ODE
    method from \pkg{OrdinaryDiffEq}.
    \item \label{enum:dense} Obtain a new dense approximation of the solution
    on time interval $[t_{k - 1}, t_k]$ by completing this discrete approximation
    using the specialized collocation polynomials in \pkg{OrdinaryDiffEq}.
    \item \label{enum:check} If in steps~\ref{enum:ode_step} and~\ref{enum:dense}
    the "old" dense approximation was only evaluated at time points $\leq t_{k - 1}$,
    we stop. Otherwise we repeat steps~\ref{enum:ode_step} and~\ref{enum:dense} iteratively based on the newly calculated approximations,
    until some stopping criterion is met.
\end{enumerate}

The fixed-point iteration can be performed efficiently since \pkg{OrdinaryDiffEq}
allows us to specify that the Jacobian and its factorization should only
be evaluated once.
Optionally, Anderson acceleration~\cite{anderson65_iterative_procedures,walker11_anderson_acceleration_fixed_point} can be used to accelerate the convergence of the fixed-point iteration.
We have observed that, e.g., for the one-dimensional IVP

\begin{equation*}
    \begin{aligned}
        x'(t) &= - x(t - 1/3) - x(t - 1/5), \qquad & t &\in [0, 100],\\
        x(t) &= \delta_{0}(t), \qquad &t &\leq 0,
    \end{aligned}
\end{equation*}

Anderson acceleration reduces the number of evaluations of the model,
fixed-point iterations, and convergence failures from $5517$, $720$, and $58$
to $2829$, $308$, and $30$, respectively, when the explicit ODE
method \alg{Tsit5}, relative tolerance $10^{-3}$, and absolute tolerance $10^{-6}$
are used.

\subsection{Discontinuities}

User-provided discontinuities and discontinuities arising from constant time delays
up to the order of the ODE method are added as grid points. Existing functionality
in \pkg{OrdinaryDiffEq} for handling discontinuities ensures that the numerical solver
steps to these discontinuities exactly and integration is restarted appropriately
at these time points.

Discontinuities due to state-dependent delays are handled by locating time points $\xi$
that satisfy
\begin{equation}\label{eq:state_dependent}
    \zeta = \xi - \tau(\tilde{x}(\xi), \xi),
\end{equation}
where $\zeta$ is a known discontinuity and $\tilde{x}$ is the current approximation of
the solution $x$, and adding them as grid points. We search for such discontinuities 
only if a step is rejected and rely on the error estimator otherwise. If a step
is rejected, we check for sign changes of $\xi_i - \tau(\tilde{x}(\xi_i), \xi_i) - \zeta_j$
at pre-defined number of equally spaced time points $\xi_i$ in the current time
interval $[t_{k - 1}, t_k]$, where $\zeta_j \leq t_{k - 1}$ are known discontinuities.
If a sign change is detected, a solution $\xi$ of eq.~\eqref{eq:state_dependent} is
calculated by a root finding algorithm and added as grid point, and the step size is
reduced accordingly. We do not try to find the earliest time point that satisfies
eq.~\eqref{eq:state_dependent} but stop as soon as one discontinuity is found in the
time interval $[t_{k - 1}, t_k]$. A very similar algorithm is used by
the DDE solver \pkg{RADAR5} by Guglielmi and Hairer~\cite{guglielmi01_implementing_radau},
however, they additionally locate discontinuities if the error estimates increases
above a certain threshold in subsequent time steps.

Users may also specify no time delays at all. In that case, the choice of a method 
with residual error estimator such as \alg{RK4} or \alg{OwrenZen5} is strongly 
recommended because the residual error estimator attempts to ensure accuracy over the
entire solution interval instead of simply at step points. This may be required for
the history function (i.e., the ODE integrator's interpolation) to have the relevant
accuracy.

\subsection{Interface}

We do not present all features and options of \pkg{DelayDiffEq} but only provide one
example for how DDEs can be declared and solved with \pkg{DelayDiffEq}. For a more
thorough explanation we refer the reader to the extensive documentation of the
\pkg{DifferentialEquations} framework.\footnote{The documentation is available
online at \url{http://docs.sciml.ai/stable/}.}

The interface of \pkg{DelayDiffEq} is similar to the one for ODEs in \pkg{OrdinaryDiffEq}
or stochastic differential equations (SDEs) in \pkg{StochasticDiffEq}~\cite{rackauckas2017adaptive}.
The code snippet below shows how a version of the so-called Hutchinson's
equation~\cite{hutchinson48_circul_causal_system_in_ecolog},
a population growth model also known as delay logistic equation, can be formulated
with \pkg{DelayDiffEq}. The model is given by
\begin{equation}\label{eq:hutchinson}
    \begin{aligned}
        x'(t) &= - x(t - 1),& \qquad t &\geq 0, \\
        x(t) &= 1,& \qquad t &\leq 0.
    \end{aligned}
\end{equation}
In \pkg{DelayDiffEq}, a DDE problem consists of the differential equation, the history function,
the integration time span, and the time delays. Optionally a value at the initial time
point and model parameters can be provided.

\begin{jllisting}
using DelayDiffEq

# Hutchinson's equation
f(x, h, p, t) = - h(p, t - 1)

# History function
h(p, t) = 1.0

# Integration time span
tspan = (0.0, 10.0)

# Complete DDE problem formulation
prob = DDEProblem(f, h, tspan;
                  constant_lags = (1.0,))
\end{jllisting}

Although not useful in this case, for illustrative purposes we demonstrate how a
problem with state-dependent delay can be formulated.

\begin{jllisting}
# Complete DDE problem formulation
# with state-dependent delay
lag(p, t) = 1.0
prob = DDEProblem(f, h, tspan;
                  dependent_lags = (lag,))
\end{jllisting}

In both cases presented above, 

\begin{jllisting}
solve(prob, MethodOfSteps(Tsit5()))
\end{jllisting}

computes a solution to the problem with the ODE method \alg{Tsit5}, using the
default solver options. \pkg{DelayDiffEq} supports DDEs with mass matrices and neutral 
DDEs as well. Similar to the ODE case, singular mass matrix DDEs generate
differential-algebraic equations (DAEs) which are then solved via the \pkg{OrdinaryDiffEq}
DAE solvers. 

From the structure of the design, \pkg{DelayDiffEq} includes all of the features of the
\pkg{OrdinaryDiffEq} ODE solvers, including (but not limited to):

\begin{itemize}
    \item Explicit Runge-Kutta methods from 1st order to 9th order, and arbitrary order extrapolation methods
    \item Semi-implicit Rosenbrock methods (allows for singular mass matrices)
    \item Implicit methods via Singularly Diagonally Implicit Runge-Kutta methods (SDIRK) and Backwards Differentiation Formulae (BDF)
    \item Fully-Implicit Runge-Kutta (FIRK) methods
    \item Stabilized explicit methods via Runge-Kutta Chebyshev methods (ROCK schemes)
    \item Dense output
    \item Support for sensitivity analysis via forward and reverse mode automatic differentiation
    \item Event handling for hybrid equations
    \item Jacobians for implicit methods defined by automatic differentiation, along with the ability to define analytical Jacobians
    \item Support for a wide range of linear solvers, such as preconditioned GMRES
\end{itemize}

We note in passing that the quasi-Newton schemes used in the ODE solvers only require the
Jacobian and its factorization at the beginning of the step for the non-FIRK implicit methods.
This means that the entire handling of the unconstrained time stepping of 
Section~\ref{sec:unconstrained} can be done without recalculating the Jacobian or refactorizing
the Jacobian. As these steps are commonly the most costly aspect for solving large stiff equations,
this greatly improves the efficiency of the implementation over using a naive Newton solve.

\section{Examples}

\subsection{Stiff problem: Antibody production}

Next we consider a DDE model of antigen stimulated antibody production by Waltman~\cite{waltman78_threshold_model_antigen}.
The model is given by
\begin{equation}\label{eq:waltman}
    \begin{aligned}
        x_1'(t) &= - r x_1(t) x_2(t) - s x_1(t) x_4(t), \\
        x_2'(t) &= - r x_1(t) x_2(t) \\
                &\qquad + \alpha r x_1(x_5(t)) x_2(x_5(t)) H(t - t_1), \\
        x_3'(t) &= r x_1(t) x_2(t), \\
        x_4'(t) &= - s x_1(t) x_4(t) - \gamma x_4(t) \\
                &\qquad +\beta r x_1(x_6(t)) x_2(x_6(t)) H(t - t_2), \\
        x_5'(t) &= H(t - t_1) \frac{x_1(t) x_2(t) + x_3(t)}{x_1(x_5(t)) x_2(x_5(t)) + x_3(x_5(t))}, \\
        x_6'(t) &= H(t - t_2) \frac{10^{-12} + x_2(t) + x_3(t)}{10^{-12} + x_2(x_6(t)) + x_3(x_6(t))},
    \end{aligned}
\end{equation}
where $H$ is the Heavyside step function, i.e., $H(t) = 0$ if $t < 0$ and
$H(t) = 1$ if $t \geq 0$. As Guglielmi and Hairer~\cite{guglielmi01_implementing_radau},
we consider the model on the time interval $t \in [0, 300]$ with parameter values
$\alpha = 1.8$, $\beta = 20$, $\gamma = 0.002$, $r = 5 \times 10^4$, $s = 10^5$,
$t_1 = 35$, $t_2 = 197$, and history functions $x_1(t) = 5 \times 10^{-6}$,
$x_2(t) = 10^{-15}$, and $x_3(t) = x_4(t) = x_5(t) = x_6(t) = 0$ for all $t \leq 0$.

\begin{figure}
    \centering
    \input{figures/waltman.tex}
    \caption{Solution of problem~\eqref{eq:waltman}.}
    \label{fig:waltman}
\end{figure}
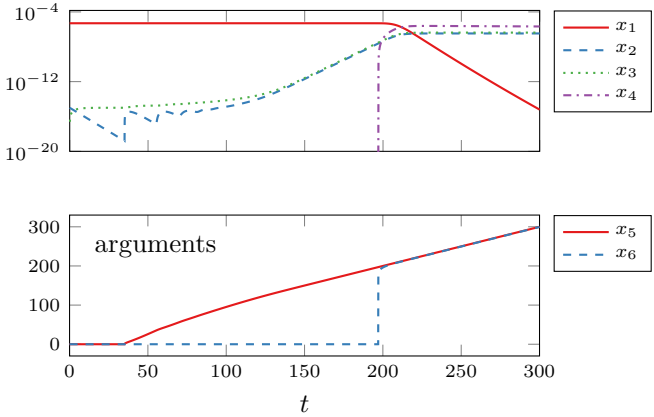

This problem is very stiff, contains vanishing state-dependent delays, and
displays steep increases in components $x_2$, $x_4$, and $x_6$ at time points
$t_1$ and $t_2$~\cite{guglielmi01_implementing_radau,waltman78_threshold_model_antigen}.
Guglielmi and Hairer~\cite{guglielmi01_implementing_radau} compared their code RADAR5
with five existing numerical solvers, and only RADAR5 was able to solve
problem~\eqref{eq:waltman}.\footnote{A later produced \href{https://www.mapleprimes.com/posts/201287-Delay-Differential-Equations--Maple-Is-A-Leader}{blog post} showed that Maple was able to successfully solve the problem, while Mathematica and MATLAB solvers were not.}
As shown in Figure~\ref{fig:waltman}, \pkg{DelayDiffEq} manages
to reproduce Figure~3 in~\cite{guglielmi01_implementing_radau} using the stiff ODE
solver \alg{KenCarp5}. As Guglielmi and Hairer~\cite{guglielmi01_implementing_radau},
we used relative tolerances of $10^{-9}$ and absolute tolerances $10^{-21}$ and
$10^{-9}$ for the first four and the last two components, respectively. We are
pleased to see that \pkg{DelayDiffEq} can cope with the discontinuities at the time points
$t_1 = 35$ and $t_2 = 197$, even though they are not declared as grid points as
in~\cite{guglielmi01_implementing_radau}.

\subsection{Sensitivity analysis via automatic differentiation}

The \pkg{SciMLSensitivity} tutorials highlight that \pkg{DelayDiffEq} is compatible with forward- and reverse-mode automatic differentiation.
For this example, we define a Lotka-Volterra model with delay and constant history function:

\jlinputlisting[linerange={1-15}]{src/sensitivity.jl}

The following function computes the solution of the DDE model on an equally-spaced time grid:

\jlinputlisting[linerange={24-36}]{src/sensitivity.jl}

This function includes the DDE solver and is compatible with forward- and reverse mode automatic differentiation:

\jlinputlisting[linerange={45-50}, xgobble=4]{src/sensitivity.jl}

\subsection{Algorithm development and benchmarking}

\begin{figure}
    \centering
    \input{figures/benchmark.tex}
    \caption{Benchmarks on a DDE model of quorum sensing~\cite{BuddrusSchiemann2014}.}
    \label{fig:benchmark}
\end{figure}
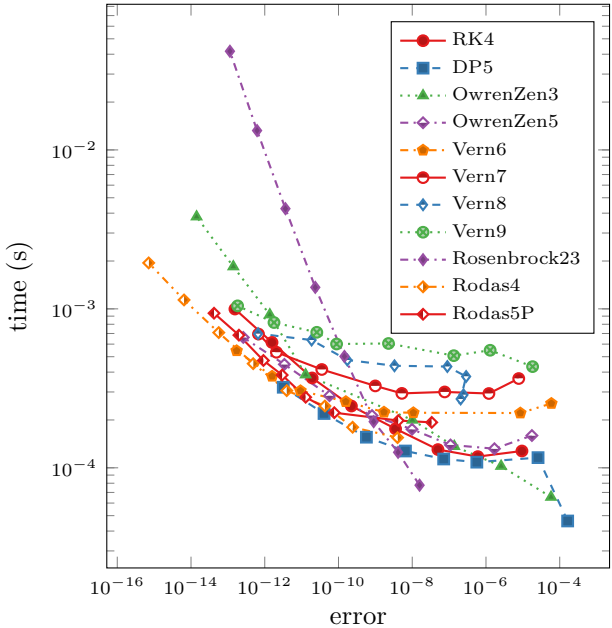

To facilitate the research and development of efficient DDE solving, \pkg{DelayDiffEq}
ties into the SciMLBenchmarks ecosystem to provide an easy way to benchmark its 
DDE solvers.
A demonstration of a work-precision diagram generated by \pkg{DiffEqDevTools}
on a delay differential equation model of quorum sensing~\cite{BuddrusSchiemann2014} is shown in Figure~\ref{fig:benchmark}.

\jlinputlisting[linerange={3-41}]{src/benchmark.jl}

Notice how with very
little code, DDE integrators for non-stiff and stiff DDEs can be compared to help
optimize the implementations.

\subsection{Stochastic delay differential equations}

This same infrastructure of \pkg{DelayDiffEq} was used on \pkg{StochasticDiffEq} for SDEs to generate convergent solvers for stochastic delay differential equations, as was demonstrated in~\cite{sykora2020stochasticdelaydiffeq}.

\section{Conclusion}

In this paper we presented \pkg{DelayDiffEq}, a Julia software package for numerically solving DDEs.
We explained the employed numerical algorithm and showed that it can successfully
solve very stiff systems with state-dependent delays. But most importantly, the design of \pkg{DelayDiffEq}
allows it to inherit all of the development work from the \pkg{OrdinaryDiffEq} solvers 
of \pkg{DifferentialEquations}, a design which in turns leads to rapid development of
\pkg{DelayDiffEq} even without developers having to focus on DDEs. This is thus an example of how
composable software designs can change the efficiency of scientific computing.

\subsection*{Acknowledgements}

Large parts of this work were completed for a Master thesis about delay differential equation models of quorum sensing at the department of Mathematics at TU Munich, Germany, under the supervision of Christina Kuttler.
This work was financially supported by the Centre for Interdisciplinary Mathematics (CIM) at Uppsala University, Sweden, and by the Swedish Foundation for Strategic Research via the project \emph{ASSEMBLE} (contract number: RIT15-0012).

This material is based upon work supported by the National Science Foundation under grant no. OAC-1835443, grant no. SII-2029670, grant no. ECCS-2029670, grant no. OAC-2103804, and grant no. PHY-2021825. We also gratefully acknowledge the U.S. Agency for International Development through Penn State for grant no. S002283-USAID. The information, data, or work presented herein was funded in part by the Advanced Research Projects Agency-Energy (ARPA-E), U.S. Department of Energy, under Award Number DE-AR0001211 and DE-AR0001222. We also gratefully acknowledge the U.S. Agency for International Development through Penn State for grant no. S002283-USAID. The views and opinions of authors expressed herein do not necessarily state or reflect those of the United States Government or any agency thereof. This material was supported by The Research Council of Norway and Equinor ASA through Research Council project "308817 - Digital wells for optimal production and drainage". Research was sponsored by the United States Air Force Research Laboratory and the United States Air Force Artificial Intelligence Accelerator and was accomplished under Cooperative Agreement Number FA8750-19-2-1000. The views and conclusions contained in this document are those of the authors and should not be interpreted as representing the official policies, either expressed or implied, of the United States Air Force or the U.S. Government. The U.S. Government is authorized to reproduce and distribute reprints for Government purposes notwithstanding any copyright notation herein.

\bibliographystyle{abbrv}
\bibliography{ref.bib}

\end{document}

%% file: figures/mackey_glass.tex
\begin{tikzpicture}
    \begin{groupplot}[
        group style={
            group size={1 by 2},
            vertical sep={5\baselineskip},
        },
        width={0.9\linewidth},
        height={0.5\linewidth},
        no marks,
    ]  
        \nextgroupplot[
            title={solution},
            xlabel={$t$},
            ylabel={$x(t)$},
            xmin={300.0},
            xmax={600.0},
        ]
        \addplot+ table [x=t,y=x] {mackey_glass.csv};
    
        \nextgroupplot[
            title={time delay embedding},
            xlabel={$x(t)$},
            ylabel={$x(t - \tau)$},
        ]
        \addplot+ table [x=x,y=xtau] {mackey_glass.csv};
    \end{groupplot}
\end{tikzpicture}

%% file: figures/waltman.tex
\begin{tikzpicture}
\begin{groupplot}[
    group style={
        group size={1 by 2},
        x descriptions at={edge bottom},
        vertical sep={2\baselineskip},
    },
    xlabel={$t$},
    xmin={0.0},
    xmax={300.0},
    width={0.9\linewidth},
    height={0.4\linewidth},
    legend pos={outer north east},
    no marks
]
    \nextgroupplot[xmode=normal,ymode=log,ymin={1.0e-20},filter discard warning={false}]
    \foreach \i in {1,...,4}{
        \addplot+ table[x=timestamp,y={value\i}] {waltman.csv};
        \addlegendentryexpanded{$x_\i$}
    }

    \nextgroupplot[]
    \foreach \i in {5,6}{
        \addplot+[] table[x=timestamp,y={value\i}] {waltman.csv};
        \addlegendentryexpanded{$x_\i$}
    }
    \node [right] at (10, 250) {arguments};
\end{groupplot}
\end{tikzpicture}

%% file: figures/benchmark.tex
\begin{tikzpicture}
\begin{axis}[
    xlabel={error},
    ylabel={time (s)},
    xmode={log},
    ymode={log},
    filter discard warning={false},
    legend pos={north east},
    width={0.95\linewidth},
    height={1.05\linewidth},
]
    \foreach \alg in {RK4,DP5,OwrenZen3,OwrenZen5,Vern6,Vern7,Vern8,Vern9,Rosenbrock23,Rodas4,Rodas5P}{
        \addplot+ table [discard if not={alg}{\alg}, x=error, y=time] {benchmark.csv};
        \addlegendentryexpanded{\alg}
    }
\end{axis}
\end{tikzpicture}

%% file: paper.bbl
\begin{thebibliography}{10}

\bibitem{anderson65_iterative_procedures}
D.~G. Anderson.
\newblock Iterative procedures for nonlinear integral equations.
\newblock {\em Journal of the {ACM}}, 12(4):547--560, 10 1965.

\bibitem{baker95_issues_numer_solut_evolut_delay_differ_equat}
C.~T.~H. Baker, C.~A.~H. Paul, and D.~R. Will{\'e}.
\newblock Issues in the numerical solution of evolutionary delay differential
  equations.
\newblock {\em Advances in Computational Mathematics}, 3(1):171--196, 1995.

\bibitem{bellen03_numer_methods_delay_differ_equat}
A.~Bellen and M.~Zennaro.
\newblock {\em Numerical Methods for Delay Differential Equations}.
\newblock Oxford University Press, 3 2003.

\bibitem{bezanson2017julia}
J.~Bezanson, A.~Edelman, S.~Karpinski, and V.~B. Shah.
\newblock Julia: A fresh approach to numerical computing.
\newblock {\em SIAM review}, 59(1):65--98, 2017.

\bibitem{BuddrusSchiemann2014}
K.~Buddrus-Schiemann, M.~Rieger, M.~M\"{u}hlbauer, M.~V. Barbarossa,
  C.~Kuttler, B.~A. Hense, M.~Rothballer, J.~Uhl, J.~R. Fonseca,
  P.~Schmitt-Kopplin, M.~Schmid, and A.~Hartmann.
\newblock Analysis of \emph{N}-acylhomoserine lactone dynamics in continuous
  cultures of \emph{Pseudomonas putida} {IsoF} by use of {ELISA} and
  {UHPLC}/{qTOF}-{MS}-derived measurements and mathematical models.
\newblock {\em Analytical and Bioanalytical Chemistry}, 406(25):6373--6383,
  Aug. 2014.

\bibitem{burrage95_parallel_sequential_methods}
K.~Burrage.
\newblock {\em Parallel and sequential methods for ordinary differential
  equations}.
\newblock Clarendon Press, Oxford, 1995.

\bibitem{dahlquist75_recent_work_stiff}
G.~Dahlquist.
\newblock Recent work on stiff differential equations.
\newblock Technical report, CM-P00069399, 1975.

\bibitem{glass79}
L.~Glass and M.~C. Mackey.
\newblock Pathological conditions resulting from instabilities in physiological
  control systems.
\newblock {\em Annals of the New York Academy of Sciences}, 316(1):214--235,
  1979.

\bibitem{guglielmi01_implementing_radau}
N.~Guglielmi and E.~Hairer.
\newblock Implementing {Radau} {IIA} methods for stiff delay differential
  equations.
\newblock {\em Computing}, 67(1):1--12, 07 2001.

\bibitem{hutchinson48_circul_causal_system_in_ecolog}
G.~E. Hutchinson.
\newblock Circular causal systems in ecology.
\newblock {\em Annals of the New York Academy of Sciences}, 50(4
  Teleological):221--246, 1948.

\bibitem{neves76_charact_jump_discont_state_depend}
K.~W. Neves and A.~Feldstein.
\newblock Characterization of jump discontinuities for state dependent delay
  differential equations.
\newblock {\em Journal of Mathematical Analysis and Applications},
  56(3):689--707, 12 1976.

\bibitem{rackauckas2017adaptive}
C.~Rackauckas and Q.~Nie.
\newblock Adaptive methods for stochastic differential equations via natural
  embeddings and rejection sampling with memory.
\newblock {\em Discrete and continuous dynamical systems. Series B},
  22(7):2731, 2017.

\bibitem{DifferentialEquations.jl-2017}
C.~Rackauckas and Q.~Nie.
\newblock {DifferentialEquations.jl} – {A} performant and feature-rich
  ecosystem for solving differential equations in julia.
\newblock {\em The Journal of Open Research Software}, 5(1), 2017.
\newblock Exported from https://app.dimensions.ai on 2019/05/05.

\bibitem{shampine01_solving_ddes_matlab}
L.~F. Shampine and S.~Thompson.
\newblock Solving {DDEs} in {Matlab}.
\newblock {\em Applied Numerical Mathematics}, 37(4):441--458, 6 2001.

\bibitem{sykora2020stochasticdelaydiffeq}
H.~T. Sykora, C.~Rackauckas, D.~Widmann, and D.~Bachrathy.
\newblock {StochasticDelayDiffeq.jl} - {A}n integrator interface for stochastic
  delay differential equations in {Julia}.
\newblock 2020.

\bibitem{walker11_anderson_acceleration_fixed_point}
H.~F. Walker and P.~Ni.
\newblock Anderson acceleration for fixed-point iterations.
\newblock {\em {SIAM} Journal on Numerical Analysis}, 49(4):1715--1735, 1 2011.

\bibitem{waltman78_threshold_model_antigen}
P.~Waltman.
\newblock A threshold model of antigen-stimulated antibody production.
\newblock {\em Theoretical Immunology}, pages 437--453, 1978.

\end{thebibliography}
